\documentclass[12pt,a4paper]{article}
\usepackage{amsmath,amssymb}
\usepackage{enumerate,a4}
\usepackage{xcolor,colortbl}
\usepackage{hyperref,color}
\usepackage{graphicx, here}
\usepackage{multicol}
\usepackage[font=footnotesize,labelfont=bf]{caption}
\usepackage{booktabs}
\usepackage{multirow}
\usepackage{tabularx}
\usepackage{siunitx}
\usepackage{tikz}

\DeclareFontFamily{OT1}{manual}{}
\DeclareFontShape{OT1}{manual}{m}{n}{ <10> manfnt }{}

\usepackage[left=2.6cm,right=2.6cm,top=2.6cm,bottom=2.6cm,,nohead,nofoot,a4paper]{geometry}

\begin{document}
	\begin{center}{\large\bf
			ADVANCING STABLE SET PROBLEM SOLUTIONS THROUGH QUANTUM ANNEALERS
		}\\[16pt]
		Janez Povh$^{\dag,\ddag}$, Dunja Pucher$^{*}$\footnote{\!\! The second author is supported by the Austrian Science Fund (FWF): DOC 78.}\\[3mm]
		{\small
			$^{\dag}$University of Ljubljana, Faculty of mechanical engineering, Slovenia \\
			$^{\ddag,}$Rudolfovo Institute, Novo Mesto, Slovenia \\
            $^{*}$University of Klagenfurt, Department of Mathematics, Austria \\
			{\tt email:} {\tt janez.povh@fs.uni-lj.si, dunja.pucher@aau.at} }\\[18pt]

	\end{center}
	
	\small{ \noindent \textbf{Abstract}\\ 
    \noindent
   We assess the performance of D-wave quantum solvers for solving the stable set problem in a graph, one of the most studied NP-hard problems. We perform computations on some instances from the literature with up to $125$ vertices and compare the quality of the obtained solutions with known optimum solutions. It turns out that the hybrid solver gives very good results, while the Quantum Processing Unit solver shows rather modest performance overall.	
	}\\[1mm]
	\\[1mm]
	{\small\noindent\textbf{Math. Subj. Class. (2020):} Primary 90C27, Secondary 81P68}

	\date{\today}

	\section{INTRODUCTION}\label{intro}

    \noindent
    Combinatorial optimization problems arise in diverse fields such as network design, logistics, scheduling, and bioinformatics. These problems involve finding the optimal solution from a finite set of possibilities and are often NP-hard; hence, unless P~=~NP, there are no polynomial time algorithms which can solve these problems exactly. Nevertheless, a considerable effort has been dedicated to the development of efficient algorithms for solving these problems.
    
    \noindent
    Exact algorithms for combinatorial optimization problems are algorithms that guarantee finding the optimal solution. However, exact solvers are often computationally expensive, and solving larger problems becomes challenging. An alternative to exact solvers are heuristics. They provide good quality solutions significantly faster than exact methods, but do not guarantee finding the optimal solution.

    \noindent
    In recent years, a new approach for solving combinatorial optimization problems has emerged based on quantum annealing, a computing paradigm that exploits the principles of quantum mechanics to solve optimization problems. Quantum annealing can be practically implemented using quantum computers. One of the leading companies in this area is D-Wave, which has developed commercial quantum annealing systems called D-Wave Quantum Processing Units (QPUs).  However, the size of problems that can be solved using D-Wave QPUs is limited by the number of qubits in their quantum processor. To handle larger problems, D-Wave also offers hybrid solvers that combine classical with quantum computing. 

    \noindent
    This work deals with the stable set problem, one of the fundamental NP-hard problems~\cite{Karp72} in combinatorial optimization. The stable set problem asks to determine the cardinality of the largest set of pairwise nonadjacent vertices in a graph. By formulating the problem as a quadratic unconstrained binary optimization (QUBO) problem, it can be solved with D-Wave's QPU and hybrid solvers. The goal of this work is to assess and compare the performance of respective solvers in terms of quality of the solutions. For this purpose, we conduct experiments on instances from the literature with up to $125$ vertices and report computational results.

    \noindent
    This paper is organized as follows. In Section~\ref{formulations} we first define the stable set problem and formulate it as a QUBO problem. In Section~\ref{related_work} we give a short overview of related work that aims to solve the stable set problem with D-Wave. The information regarding our implementation on D-Wave, considered instances as well as computational results are presented in Section~\ref{experiments}. Finally, we conclude with a short discussion of our results and possible future work in Section~\ref{conclusions}.

    \section{THE STABLE SET PROBLEM}\label{formulations}

    \noindent
    Let $G = (V, E)$ be a simple undirected graph with $\vert V \vert = n$ vertices and $\vert E \vert = m$ edges. A set $S \subseteq V$ is called a stable set if the vertices in $S$ are pairwise nonadjacent. The stability number of $G$ is the cardinality of the largest stable set in $G$ and is denoted by $\alpha(G)$. The stable set problem asks to determine $\alpha(G)$. 

    \noindent
    A standard method in combinatorial optimization is to write an NP-hard problem as an integer program. Let $x_i$ be a binary variable indicating whether the vertex $i$ is in the stable set or not. Then the stability number of a graph $G$ is the optimal value of the following integer optimization problem with linear constraints
    \begin{equation}\label{stable_set_linear}
    \begin{aligned}
    \alpha(G) ~=~ \max \{e^Tx \mid x_i + x_j \leq 1 ~\forall \{i,j\} \in E, ~x_i \in \{0, 1\} ~\forall i \in V \}, 
    \end{aligned}
    \end{equation}
    where $e$ denotes the all-ones vector. Linear constraints in the formulation~(\ref{stable_set_linear}) can be replaced by quadratic ones. Let $x_i$ be a binary variable as defined in~(\ref{stable_set_linear}). Then the edge constraints can be written as $x_ix_j = 0$ for all $\{i,j\} \in E$, and therefore
    \begin{equation}\label{stable_set_quadratic}
    \begin{aligned}
    \alpha(G) ~=~ \max \{e^Tx \mid x_ix_j = 0 ~\forall \{i,j\} \in E, ~x_i \in \{0, 1\} ~\forall i \in V \}. 
    \end{aligned}
    \end{equation}

    \noindent
    The formulations~(\ref{stable_set_linear}) and (\ref{stable_set_quadratic}) have $n$ binary variables and $m$ constraints. The number of constraints can be reduced in the following way. Let $A$ be the adjacency matrix of a graph $G$. Since $x$ is a binary variable we have that $x_i^2 = x_i$, so the edge constraints $x_ix_j = 0$ can be written as $\frac{1}{2}x^TAx = 0$. Hence, 
    \begin{equation}\label{stable_set_quadratic_2}
    \begin{aligned}
    \alpha(G) ~=~ \max \{e^Tx \mid \frac{1}{2}x^TAx = 0, ~x_i \in \{0, 1\} ~\forall i \in V \}. 
    \end{aligned}
    \end{equation}
    
    \noindent
    In order to use D-Wave solvers, we need a QUBO formulation for the stable set problem. Generally, a QUBO is a problem of finding
    \begin{equation*}
    \begin{aligned}
    z^* ~=~ \min \{f^Tx + x^TQx \mid x \in \{0, 1\}^n \}. 
    \end{aligned}
    \end{equation*}
    
    \noindent
    In the formulation~(\ref{stable_set_quadratic_2}) we use that $\max e^Tx = - \min -e^Tx$. Furthermore, we have that the edge constraints can be written as $\frac{1}{2}x^TAx = 0$. Nevertheless, the penalty for an edge should be greater than the contribution of a vertex. That is why we impose the penalties as $\beta x^TAx$ for $\beta \geq 1$, and get to the following formulation
    \begin{equation}\label{QUBO_1}
    \begin{aligned}
    \alpha(G) ~=~ \min \{-e^Tx + \beta x^TAx \mid x \in \{0, 1\}^n \}. 
    \end{aligned}
    \end{equation}

    \noindent
    Analogous formulations with different values of the parameter $\beta$ have already been introduced. For instance, in~\cite{GHGG} the parameter was set as $\beta = 1$, while $\beta = \frac{1}{2}$ was considered in~\cite{Abello}. An alternative approach is to reward the contribution of vertices instead of penalizing the edges. Such formulation was given in~\cite{Parekh}.

    \noindent
    D-Wave solvers are designed to minimize a QUBO formulated as
    \begin{equation*}\label{QUBO_2}
    \begin{aligned}
    \min \quad & x^TQx,
    \end{aligned}
    \end{equation*}
    where $x$ is a binary variable. Due to the fact that $x_i^2 = x_i$, we can write $-e^Tx$ as $x^T(-I)x$, where $I$ denotes the identity matrix. Thus, we define for $\beta \geq 1$ the input matrix $Q$ as
    \begin{equation}\label{matrix_Q}
    \begin{aligned}
    Q = -I + \beta A.
    \end{aligned}
    \end{equation}

    \noindent
    It is important to note that the solution obtained from D-Wave is not necessary the optimal solution. But also, we would like to emphasize that the solution returned by D-Wave is not necessarily a stable set, either. Although we impose certain penalties in~(\ref{QUBO_1}), it does not mean that the sum of penalties in the final solution will be zero. The correctness of the solution can be easily checked. If
    \begin{equation*}
    \begin{aligned}
    \vert e^Tx \vert \neq \vert -e^Tx + \beta x^TAx \vert,
    \end{aligned}
    \end{equation*}
    then the underlying solution is clearly not a stable set. We will discuss this issue in Section~\ref{experiments} in more detail.

    \section{RELATED WORK}\label{related_work}

    The previous work focused on solving the stable set problem with the D-Wave QPU solver. For this purpose, several authors considered instances which fit the architecture of D-Wave, while for other arbitrary instances partitioning methods were proposed. The performance was compared with exact methods or heuristics. The authors focused mostly on computational times and investigated whether a speed-up can be realized or not. In the following text we give more details regarding the respective research.

    \noindent
    One of the first practical results for solving the stable set problem with quantum annealing in D-Wave was shown in~\cite{Parekh}. The results were obtained on D-Wave Two AQO with $512$ qubits in Chimera graph architecture. The experiments were done on Chimera graphs, which fit the underlying architecture. However, it turned out that Selby's exact method~\cite{selby:2014} outerperformed the quantum annealing approach. 

    \noindent
    Further numerical experiments were conducted in~\cite{GHGG}, and that on D-Wave 2X with about $1000$ qubits in the Chimera graph architecture. The first set of experiments was done on random graphs with $45$ vertices. It turned out that classical solvers available are usually faster for such small instances. The second set of experiments was performed with subgraphs of D-wave chimera graphs. For large instances with $800$ or more vertices, D-Wave provided the best solutions. The last set of experiments was performed with instances that generally do not fit the D-Wave architecture. A decomposition method was proposed for these instances. However, it turned out that the proposed method was effective only for relatively thin graphs. The method in question was further improved in\cite{Pelofske:2019}.
    
    \noindent
    Experiments on D-Wave 2000Q were done in~\cite{YPB}. The obtained results for random graphs with up to $60$ vertices were compared to some classical algorithms. It was shown that D-Wave QPU was outerperformed by simulated thermal annealing.

    \noindent
    Recently, the algorithm introduced in~\cite{Pelofske:2019} was combined with the method of parallel quantum annealing~\cite{Pelofske:2022xho}. The suggested approach was applied on graphs with up to $120$ vertices. All computations were performed on the D-Wave Advantage System 4.1. Among others, it was shown that the proposed algorithm can compute solutions for certain instances up to around two to three orders of magnitude faster than a classical solver. 
     
    \section{COMPUTATIONAL RESULTS}\label{experiments}

   In this section we report computational experiences and compare the solutions obtained by D-Wave with the exact ones, known from the literature or recomputed by BiqBin solver \cite{Gus:2022}. Note that all codes and data are available as supplementary files on the arXiv page of this paper and can also be found on github\footnote{\url{https://github.com/DunjaPucher/Stable_set_experiments.}}.

    \noindent
    For our computations we consider instances from the literature with up to $125$ vertices. More precisely, we consider some instances from the Second DIMACS Implementation Challenge~\cite{Johnson1996Cliques}. Furthermore, we consider some small Paley graphs, torus graphs as well as spin graphs. The stability numbers of these instances were previously considered in~\cite{Gaar:2022}. 

    \noindent
    We perform computations on Leap quantum cloud service by using the class of Binary Constraint Models (BQM), since the problems in this class can be solved with both the QPU as well as the hybrid solver. We use D-Wave's Advantage System, which uses the Pegasus topology and contains more than $5000$ qubits~\cite{boothby:2020}. Finally, we would like to note that we perform all computations with default parameters, and for QPU computations we set the number of reads to be $1000$.

    \noindent
    The computational results obtained by QPU and hybrid solvers are gathered in Tables~\ref{Table_results_QPU} and~\ref{Table_results_hybrid}. The columns $1-4$ contain general information about the instances: name of the instance, number of vertices $n$, number of edges $m$, as well as the stability number $\alpha(G)$. For each instance we perform several experiments for different values of the parameter $\beta$ in~(\ref{matrix_Q}). We set $\beta = 1$, $\beta = 10$ and $\beta = 100$. The columns $5-7$ contain the information about the cardinalities of the solutions obtained from the solver, i.e.\ about the stability numbers of the considered instance. As previously mentioned, the solution returned by the solver is not necessarily a stable set. Therefore, we check whether the solution is a stable set or not. If the solver returned a solution which is not a stable set, we mark the cell as dark gray. If the solver returned a solution which is not optimal, we mark the cell as gray. 

\begin{table}
\footnotesize
\begin{center}
\caption{Results for computations done with QPU solver}
\begin{tabular}{ l r r r r r r r}
Graph & $n$ & $m$ & $\alpha(G)$ & $\beta = 1$ & $\beta = 10$ & $\beta = 100$  \\
\hline 
  C125.9 & 125 & 787 & 34 & \cellcolor{lightgray}24 & \cellcolor{lightgray}22 & \cellcolor{lightgray}20 \\
 \hline
 DSJC125.5 & 125 & 3859 & 10 & \cellcolor{lightgray}4 & \cellcolor{lightgray}4 & \cellcolor{gray}6 \\
 DSJC125.9 & 125 & 789 & 34 & \cellcolor{gray}26 & \cellcolor{lightgray}21 & \cellcolor{lightgray}19 \\
 \hline
 hamming6\_2 & 64 & 192 & 32 & 32 & \cellcolor{lightgray}29 & \cellcolor{lightgray}25 \\
 hamming6\_4 & 64 & 1312& 4 & 4 & \cellcolor{lightgray}3 & 4 \\
 \hline
 johnson8\_2\_4 & 28 & 168 & 4 & 4 & 4 & 4 \\
 johnson8\_4\_4 & 70 & 560 & 14 & \cellcolor{lightgray}11 & \cellcolor{lightgray}9 & \cellcolor{lightgray}9 \\
 johnson16\_2\_4 & 120 & 1680 & 8 & \cellcolor{gray}8 &  \cellcolor{lightgray}6 & \cellcolor{lightgray}7 \\
 \hline
 MANN\_a9 & 45 & 72 & 16 & 16 & \cellcolor{lightgray}15 & \cellcolor{lightgray}15 \\
 \hline
 paley61 & 61 & 915 & 5 & 5 & 5 & 5  \\
 paley73 & 73 & 1314 & 5 & 5 & 5 & \cellcolor{lightgray}4 \\
 paley89 & 89 & 1958 & 5 & \cellcolor{lightgray}4 & 5 & 5\\
 paley97 & 97 & 2328 & 6 & \cellcolor{lightgray}5 & \cellcolor{lightgray}5 & \cellcolor{lightgray}4 \\
 paley101 & 101 & 2525 & 5 & 5 & 5 & 5 \\
 \hline
 spin5 & 125 & 375 & 50 & \cellcolor{lightgray}46 & \cellcolor{lightgray}35 & \cellcolor{lightgray}34 \\
 \hline
 torus11 & 121 & 242 & 55 & 55 & \cellcolor{lightgray}43 & \cellcolor{lightgray}40 \\
 \hline
\end{tabular}
\label{Table_results_QPU}
\end{center}
\end{table}

\begin{table}
\begin{center}
\caption{Results for computations done with hybrid solver}
\footnotesize
\begin{tabular}{ l r r r r r r r}
Graph & $n$ & $m$ & $\alpha(G)$ & $\beta = 1$ & $\beta = 10$ & $\beta = 100$  \\
\hline 
  C125.9 & 125 & 787 & 34 & 34 & 34 & 34 \\
 \hline
 DSJC125.5 & 125 & 3859  & 10 & 10 & 10 & 10 \\
 DSJC125.9 & 125 & 789 & 34 & 34 & 34 & 34 \\
 \hline
 hamming6\_2 & 64 & 192 & 32 & 32 & 32 & 32 \\
 hamming6\_4 & 64 & 1312 & 4 & 4 & 4 & 4 \\
 \hline
 johnson8\_2\_4 & 28 & 168 & 4 & 4 & 4 & 4 \\
 johnson8\_4\_4 & 70 & 560 & 14 & 14 & 14 & 14 \\
 johnson16\_2\_4 & 120 & 1680 & 8 & 8 & 8 & 8 \\
 \hline
 MANN\_a9 & 45 & 72 & 16 & 16 & 16 & 16 \\
 \hline
 paley61 & 61 & 915 & 5 & 5 & 5 & 5 \\
 paley73 & 73 & 1314 & 5 & 5 & 5 & 5 \\
 paley89 & 89 & 1958 & 5 & 5 & 5 & 5 \\
 paley97 & 97 & 2328 & 6 & 6 & 6 & 6\\
 paley101 & 101 & 2525 & 5 & 5 & 5 & 5 \\
 \hline
 spin5 & 125 & 375 & 50 & 50 & 50 & 50 \\
 \hline
 torus11 & 121 & 242 & 55 & 55 & 55 & 55 \\
 \hline
\end{tabular}
\label{Table_results_hybrid}
\end{center}
\end{table}

\noindent
Based on the results in Table~\ref{Table_results_QPU} we make several observations. First, not every solution returned by QPU is a stable set, see results for \textit{dsjc125.5}, \textit{dsjc125.9} and \textit{johnson16\_2\_4}. Second, obtained solutions are not always optimal; some of the solutions are quite far from the optimal values, see for example results for \textit{C125.9} and \textit{dsjc} graphs. The smallest instance for which we got a non-optimal value for parameter $\beta = 10$ is \textit{hamming6\_4} with $64$ vertices. Finally, we observe that imposing higher penalties on edges, i.e.\ setting the value of the parameter $\beta$ to $10$ or $100$ does not necessarily yield better solutions than setting $\beta = 1$, see for instance results for \textit{dsjc125.9}, \textit{spin5} and \textit{torus11}.

\noindent
On the other side, from the results in Table~\ref{Table_results_hybrid} we note that for all instances and for all considered values of the parameter $\beta$ we obtained optimal solutions with the hybrid solver. 

\noindent
We now examine the solutions obtained from the QPU solver in more detail. For that purpose, we check for all gathered samples whether the solution is a stable set or not. As previously mentioned, we perform $1000$ runs per instance. Table~\ref{Feasible_solutions} contains information about the percentage of samples which are stable sets, i.e.\ solutions obtained from the QPU solver for different values of the parameter $\beta$.

\begin{table}
\begin{center}
\caption{Percentage of samples which are stable sets obtained from QPU solver}
\footnotesize
\begin{tabular}{ l r r r | l r r r }
\\
Graph & $\beta = 1$ & $\beta = 10$ & $\beta = 100$ &  Graph & $\beta = 1$ & $\beta = 10$ & $\beta = 100$ \\
\hline
C125.9 & 2.80 & 20.10 & 3.30 & MANN\_a9 & 85.31 & 95.80 & 97.20 \\
DSJC125.5 & 3.30 & 0.80 & 0.00 & paley61 & 2.00 & 7.40 & 8.70 \\
DSJC125.9 & 2.50 & 17.20 & 19.30 & paley73 & 2.20 & 5.20 & 1.70 \\
hamming6\_2 & 99.77 & 90.50 & 94.90 & paley89 & 1.30 & 3.70 & 4.90 \\
hamming6\_4 & 2.30 & 8.80 & 12.40 & paley97 & 1.70 & 5.10 & 1.60 \\
johnson8\_2\_4 & 30.90 & 49.94 & 60.96 & paley101 & 4.10 & 3.70 & 3.00 \\
johnson8\_4\_4 & 4.10 & 16.30 & 19.90 & spin5 & 32.70 & 68.00 & 61.60 \\
johnson16\_2\_4 & 0.50 & 4.30 & 3.80 & torus11 & 85.60 & 94.10 & 96.00 \\
\hline
\end{tabular}
\label{Feasible_solutions}
\end{center}
\end{table}

\noindent
The results presented in Table~\ref{Feasible_solutions} show that, for some instances and in that for all considered values of the parameter $\beta$, the percentage of samples which are stable sets is quite low. For instance, for $\beta = 1$ the percentage is under $5\%$ for $11$ out of $16$ considered instances. The situation is slightly better for $\beta = 10$ and $\beta = 100$, where this is the case for $4$, i.e.\ $7$ instances. Finally, we note again that imposing higher values of $\beta$ does not guarantee a better quality of solutions, see for instance results for \textit{dsjc125.5} and \textit{hamming6\_2}.

\section{CONCLUSIONS AND FUTURE WORK}\label{conclusions}

The goal of this work is to analyze the performance of D-Wave solvers for solving the stable set problem in a graph. We therefore used the class of BQM problems, for which a QUBO model of the problem is used. Since a stable set is a subset of vertices which are pairwise not adjacent, a QUBO formulation for the stable set problem requires imposing certain penalties on the edges. We considered different levels of penalties and analyzed the quality of the solutions returned by QPU and hybrid solvers.

\noindent
We note that a general drawback of using a QUBO formulation is that basically every binary vector is a feasible solution. Hence, the obtained solution may not be a stable set at all. This fact was observed for several solutions returned by the QPU solver. Also, it appears that imposing higher penalties on the edges does not necessarily improve the quality of the solutions obtained with the QPU solver. Overall, we note that the QPU solver shows a rather modest performance, since the returned solutions are in some cases quite far from the optimum. On the other side, all results obtained by the hybrid solver were optimal. 

\noindent
In this work we focused on the instances which can be solved with both the QPU and the hybrid solver. Since the hybrid solver showed very good results, it would be interesting to investigate in detail how it performs on larger instances from the literature.

\bibliographystyle{plain}
\bibliography{References}

\end{document}